\documentclass{article}

\usepackage{graphicx}
\usepackage{amsmath,amssymb}
\usepackage{tikz}
\usepackage{eqnarray}
\usepackage{xcolor}

\usepackage{hyperref}
\newtheorem{theorem}{Theorem}[section]
\newtheorem{lemma}{Lemma}[section]
\newtheorem{proposition}{Proposition}[section]
\newtheorem{corollary}{Corollary}[section]
\newtheorem{example}{Example}[section]
\newtheorem{definition}{Definition}[section]
\newtheorem{algorithm}{Algorithm}
\usepackage{latexsym}

\def\Rmax{\mathbb{R}_{\max}}

\def\R{\mathbb{R}}

\def\Rp{\R_+}
\def\Rpn{\Rp^n}

\def\TP{\operatorname{TP}}

\def\bunity{\mathbf{1}}

\def\bzero{\mathbf{0}}

\def\opt{\operatorname{opt}} 
\def\Pol{\operatorname{P}}
\def\TPol{\operatorname{TP}}

\begin{document}
\title{Tropical analogues of a Dempe-Franke bilevel optimization problem  }
%
%
\author{Serge\u{\i} Sergeev\thanks{Supported by EPSRC grant EP/P019676/1}
and Zhengliang Liu\thanks{Emails of authors: sergiej@gmail.com, zliu082@gmail.com} \\
{\em University of Birmingham, School of Mathematics, Edgbaston B15 2TT, UK}\\
}
\date{}

\maketitle              

\begin{abstract}
We consider the tropical analogues of a particular bilevel 
optimization problem studied by Dempe and Franke~\cite{DF} and suggest 
some methods of solving these new tropical bilevel optimization problems. 
In particular, it is found that 
the algorithm developed by Dempe and Franke can be formulated and its validity can be proved in a more general setting, which includes the tropical bilevel optimization problems in question.  We also show how the feasible set can be decomposed into a finite number of tropical polyhedra, to which the tropical linear programming solvers can be applied. 

\end{abstract}

\section{Introduction}

Bilevel programming problems are hierarchical optimization problems with two levels, each of which is an optimization problem itself. 
The upper level problem models the leader's decision making problem whereas the lower level problem models the follower's problem. These two problems are coupled through common variables. 

Consider a particular problem formulated by Dempe and Franke~\cite{DF}:
\begin{equation}
\label{e:DFproblem}
\begin{split}
\min_{x,y}\quad & a^Tx+b^Ty\\
\text{s.t.}\quad &x\in \Pol_1,\ y\in  \Pol_2,\\
&y\in \arg\min_{y'}\{x^Ty'\colon y'\in\Pol_2\}.
\end{split}
\end{equation}
Here $\Pol_1$ and $\Pol_2$ are polyhedra in $\R^n$, commonly given as solution sets to some systems of affine inequalities. 

Our goal is to study some analogues and generalisations of problem~\eqref{e:DFproblem} over the tropical (max-plus) semiring. This is a special case of a general idea to develop tropical bilevel optimization, inspired both by the well-developed methodology of bilevel optimization and some of the recent successes in tropical convexity and tropical optimization~\cite{AllSimplex,AllDoubleD,GKS}.  

The tropical semiring $\Rmax=(\mathbb{R}\cup\{-\infty\}, \oplus, \otimes)$ 
is the set of real numbers $\R$ with $-\infty$, equipped  with the ``tropical addition" $\oplus$, which is taking the maximum of two numbers, and ``tropical multiplication", which is the ordinary addition~\cite{But}. Thus we have: $a \oplus b := \max(a,b)$ and $a \otimes b := a+b$, and the elements $\bzero:=-\infty$, respectively $\bunity:=0$, are neutral with respect to $\oplus$ and $\otimes$. These arithmetical operations are then extended to matrices and vectors in the usual way, and the $\otimes$ sign for multiplication will be consistently omitted. Observe that we have $a\geq\bzero$ for all $a\in\Rmax$ and hence, for example, if $c\leq d$ for some 
$c,d\in\Rmax^n$, then we have $c^Tx\leq d^Tx$ for the {\bf tropical scalar products} of these vectors with any $x\in\Rmax^n$ (note that $c^Tx$ now means 
$\max_{i=1}^n c_i+x_i$ and all matrix-vector products are understood tropically).


The maximization and minimization problems are not equivalent in 
tropical mathematics. This is intuitively clear since only one of these 
operations plays the role of addition and the other is ``dual" to it. Namely,
the maximization problems are usually easier.
Therefore, the following {\bf four} problems can be all considered as tropical analogues 
of~\eqref{e:DFproblem}.\\[1.2 ex]
{\bf  Min-min problem} (or) {\bf Max-min problem:}
\begin{equation*}
\begin{split}
&\min_{x,y}\  a^Tx\oplus b^Ty \quad\text{(or)}\quad \max_{x,y}\  a^Tx\oplus b^Ty\\
&\text{s.t.}\quad x\in \TPol_1,\  y\in \arg\min_{y'}\{x^Ty'\colon y'\in\TPol_2\},
\end{split} 
\end{equation*}
%
{\bf Min-max problem} (or) {\bf Max-max problem:}
\begin{equation*}
\begin{split}
&\min_{x,y}\  a^Tx\oplus b^Ty\quad\text{(or)}\quad  \max_{x,y}\  a^Tx\oplus b^Ty\\
&\text{s.t.}\quad x\in \TPol_1,\ y\in \arg\max_{y'}\{x^Ty'\colon y'\in\TPol_2\},
\end{split} 
\end{equation*}
where $a$ and $b$ are vectors with entries in $\Rmax$ and $\TPol_1$ and $\TPol_2$ are tropical polyhedra of $\Rmax^n$, in the sense 
of the following definition.

\begin{definition}[Tropical Polyhedra and Tropical Halfspaces]
{\rm Tropical polyhedron} is defined as an intersection of finitely many {\rm tropical affine halfspaces} defined as 
$$\{x\in\Rmax^n\mid a^Tx\oplus \alpha \leq b^T x\oplus\beta \},$$
for some $\alpha,\beta\in\Rmax$ and $a,b\in\Rmax^n$.
\end{definition}

Note that unlike the classical halfspace, the tropical halfspace is defined as a solution set of a two-sided inequality, since we cannot move terms in the absence of (immediately defined) tropical subtraction. Also note that any tropical 
polyhedron can be defined as a set of the form
\begin{equation*}
\{x\in\Rmax^n\mid  Ax\oplus c \leq Bx\oplus d \}
\end{equation*}
where $A,B$ are matrices and $c,d$ are vectors with entries in $\Rmax$ of 
appropriate dimensions. Furthermore, any tropical polyhedron is a tropically convex set in the sense of the following definition:

\begin{definition}[Tropical Convex Set and Tropical Convex Hull]
A set $C\subseteq\Rmax^n$ is called {\rm tropically convex} if for any two points $x,y\in C$, $\lambda\oplus\mu=\bunity$ then $\lambda x\oplus\mu y\in C$.

$C$ is called the {\em tropical convex hull} of $X$ if any point of $C$ is a  tropical convex combination of the points of $X$.
\end{definition}

Furthermore, it is well-known that any compact tropical polyhedron 
$C\subseteq\Rmax^n$ is the tropical convex hull of a finite number of points (e.g., \cite{AllDoubleD}).

\section{The min-min and max-min  problem}
The direct analogue of Problem \ref{e:DFproblem} is the min-min problem, which 
we consider together with the max-min problem. Here and below, the notation 
``$\opt$'' will stand for maximization or minimization. Instead of the performance measure $a^Tx\oplus b^Ty$ we will
consider a more general function $f(\cdot,\cdot)\colon \Rmax^n\times\Rmax^n\mapsto \Rmax$, for which certain properties will be assumed, depending on the situation.

Thus we consider the following problem:
\begin{equation}
\label{e:DFtrop}
\begin{split}
\opt_{x,y}\quad &f(x,y)\\
\text{s.t.}\quad &x\in \TPol_1,\ y\in \arg\min_{y'}\{x^Ty'\colon y'\in\TPol_2\},
\end{split}
\end{equation}

Using 
$\phi(x):=\min\limits_y\{x^Ty\colon y\in\TPol_2\}$
we can rewrite the lower level value function (LLVF) reformulation of~\eqref{e:DFtrop}:
\begin{equation}
\label{e:LLVF}
\begin{split}
\opt_{x,y}\quad &f(x,y)\\
\text{s.t.}\quad &x\in\TPol_1,\ y\in\TPol_2,\ x^Ty\leq \phi(x).
\end{split}
\end{equation}

Further we will assume that $f(x,y)$ is continuous and $\TPol_1$ and $\TPol_2$ are compact in the topology\footnote{In other words, $e^{f(x,y)}$ is continuous and the sets
$\{y\in\Rpn \colon \log(y)\in \TPol_1\}$ and $\{z\in\Rpn \colon \log(z)\in \TPol_1\}$ are compact in the usual Euclidean topology.} induced by the metric 
$\rho(x,y)=\max_i |e^{x_i}-e^{y_i}|$.

Let us now introduce the following notion.

\begin{definition}[Min-Essential Sets]
\label{d:miness}
Let $\TPol$ be a tropical polyhedron. Set $\mathcal{S}$ is called
a {\rm min-essential} subset of $\TPol$, if for any $x\in\Rpn$ the minimum
$\min_z\{x^Tz\colon z\in\TPol\}$ is attained at a point of $\mathcal{S}$.
\end{definition}

\begin{lemma}
\label{l:inclusion}
If $\mathcal{S}\subseteq \TPol$ is a min-essential set of $\TPol$ and $\mathcal{S}_1\subseteq\mathcal{S}_2\subseteq\TPol$ then $\mathcal{S}_2$ is also 
min-essential.
\end{lemma}

Inspired by Dempe and Franke~\cite{DF} we suggest to generalize their algorithm in order to solve~\eqref{e:DFtrop}in the form of~\eqref{e:LLVF}. Here $\mathcal{S}_{\min}(\TPol_2)$ denotes a min-essential subset of $\TPol_2$.

\begin{algorithm}
\label{a:DF}
{\bf (Solving Min-min Problem and Max-min Problem)}\\[1.2 ex]
{\rm 1. {\bf Initial step.} Find a pair $(x^0,y^0)$ solving the relaxed problem
\begin{equation}
\label{e:DFrelaxed0}
\begin{split}
&\opt_{x,y} \quad f(x,y)\\
&\text{s.t.}\quad x\in\TPol_1,\quad y\in\TPol_2.
\end{split}
\end{equation}
We verify whether $y^0\in\arg\min_{y'}\{x_0^Ty'\colon y'\in\TPol_2\}$. If ``yes" then {\bf stop}, 
$(x^0,y^0)$ is a solution.\\
If not then find a point $z^0$ of $\mathcal{S}_{\min}(\TPol_2)$ that attains $\min_{y'}\{(x^0)^Ty'\colon y'\in\TPol_2\}$.
Let $\mathcal{Z}^{(0)}=\{z^0\}$.\\[1.2 ex]
2. {\bf General step.}  Find a pair $(x^k,y^k)$ solving the problem
\begin{equation}
\label{e:DFrelaxed}
\begin{split}
&\opt_{x,y} f(x,y)\\
& \text{s.t.}\quad x\in\TPol_1,\quad y\in\TPol_2,\quad x^Ty\leq \min\limits_{z\in \mathcal{Z}^{(k-1)}} x^Tz.
\end{split}
\end{equation}
We verify whether $y^k\in\arg\min_{y'}\{(x^k)^Ty'\colon y'\in\TPol_2\}$. If ``yes" then {\bf stop}, 
$(x^k,y^k)$ is a solution.\\
If not then find a point $z^k\in \mathcal{S}_{\min}(\TPol_2)$ that attains $\min_{y'}\{(x^k)^Ty'\colon y'\in\TPol_2\}$.
Let $\mathcal{Z}^{(k)}=\mathcal{Z}^{(k-1)}\cup \{z^{k}\}$, and repeat 2. with $k:=k+1$.} $\hfill\square$
\end{algorithm}

We now include the proof of convergence and validity of this algorithm, 
although it just generalizes the one given by Dempe and Franke~\cite{DF}.

\begin{theorem}
\label{p:DF}
Let $\mathcal{S}_{\min}(\TPol_2)$ be finite. 
Then Algorithm~\ref{a:DF} terminates in a finite number of steps
and results in a globally optimal solution
of~\eqref{e:DFproblem}.     
\end{theorem}
\emph{Proof}
First observe that as $\TPol_1$ and $\TPol_2$ are compact then
the feasible set of~\eqref{e:DFrelaxed0} is also compact. 
The feasible set of~\eqref{e:DFrelaxed} is also compact as 
intersection of the compact set $\TPol_1\times\TPol_2$ with the 
closed set
\begin{equation}
\label{e:newconstraints}
\{(x,y)\colon x^Ty\leq x^Tz\;\quad \forall z\in \mathcal{Z}^{(k-1)}\}.
\end{equation}
As $f(x,y)$ is continuous as a function of $(x,y)$,
the optima in~\eqref{e:DFrelaxed0} and~\eqref{e:DFrelaxed} always exist.

Now consider the sequence $\{z^k\}_{k=0}^{\infty}$ generated by 
the algorithm. 
Points $z^k$ belong to a 
finite (min-essential) subset of $\TPol_2$ and hence 
there exist $k_1$ and $k_2$ such that $k_1<k_2$ and $z^{k_1}=z^{k_2}$.
However, $z^{k_1}\in\mathcal{Z}^{(k_2-1)}$ and hence
$$
\min\limits_{z\in\mathcal{Z}^{(k_2-1)}} (x^{k_2})^T z\leq 
(x^{k_2})^Tz^{k_1}=\min\limits_{z\in\TP_2} (x^{k_2})^T z\leq 
\min\limits_{z\in\mathcal{Z}^{(k_2-1)}} (x^{k_2})^T z.
$$
The inequalities turn into equalities, and $(x^{k_2},z^{k_2})$ is a globally 
optimal solution since it is feasible for~\eqref{e:DFtrop} and globally optimal 
for its relaxation~\eqref{e:DFrelaxed}. $\hfill\blacksquare$

Let us now argue that a finite min-essential set exists for each
tropical polyhedron $\TPol$.

\begin{definition}[Minimal Points]
\label{d:minimal}
Let $\TPol$ be a tropical polyhedron.
A point $x\in\TPol$ is called {\rm minimal} if 
$y\leq x$ and $y\in\TPol$ imply $y=x$. 
The set of all minimal points of $\TPol$
is denoted by $\mathcal{M}(\TPol)$.
\end{definition}

\begin{definition}[Extreme Points]
\label{d:extremal}
Let $\TP$ be a tropical polyhedron.
A point $x\in\TP$ is called {\rm extreme} if any equality
$x=\lambda u\oplus\mu v$ with $\lambda\oplus\mu=\bunity$ and 
$u,v\in\TP$ implies $x=u$ or $x=v$. 
\end{definition}

We have the following known observation. Note, however, that this observation does not hold in the usual convexity, as counterexamples on the plane can be easily constructed.  

\begin{lemma}[Helbig~\cite{Helbig}]
\label{l:minimal-extremal}
Any minimal point of a tropical polyhedron is extreme.
\end{lemma}

The set of extreme points of a tropical polyhedron is finite, see for example 
Allamigeon, Gaubert and Goubault~\cite{AllDoubleD}. Combining this 
with an observation that the set 
$\{z\in\TPol\colon z\leq y\}$ is compact 
and hence contains a minimal point, we obtain the following claims.

\begin{proposition}
\label{p:minimal-fes}
$\mathcal{M}(\TPol)$ is a finite (and non-empty) 
min-essential subset for any tropical polyhedron $\TPol$. 
\end{proposition}

\begin{corollary}
\label{c:ess-exists}
Any tropical polyhedron has a finite min-essential subset.
\end{corollary}

Several problems arise when trying to implement the general 
Dempe-Franke algorithm in tropical setting. One of them is how to find a point of a finite min-essential set $\mathcal{S}_{min}(\TPol_2)$ that attains $\min_{y'}\{(x^k)^Ty'\colon y'\in\TPol_2\}$ and which min-essential set to 
choose. An option here is to exploit 
the tropical simplex method of Allamigeon, Benchimol, Gaubert and Joswig~\cite{AllSimplex}, which (under some generically true conditions imposed on $\TPol_2$) can find a point that attains
$\min_{y'}\{(x^k)^Ty'\colon y'\in\TPol_2\}$ 
and belongs to the set of tropical basic points of 
$\TPol_2$. The set of tropical basic points is finite and includes all extreme points~\cite{AllSimplex} and hence all the minimal points of $\TPol_2$, thus it is also a finite min-essential subset of $\TPol_2$ by Lemma~\ref{l:inclusion}.

Even more imminent problem is how to solve~\eqref{e:DFrelaxed}, as 
the techniques referred to in Dempe and Franke~\cite{DF} are not immediately "tropicalized". An option here is to use reduction of the constraints defining a tropical polyhedron to MILP constraints. Such reduction was suggested, e.g., in De Schutter, Heemels and Bemporad~\cite{Bart3} based on \cite{Bart2}. 
More precisely, we need to consider constraints of the following two kinds: 1) $a^Tx\leq  \alpha$ and 2) $a^Tx\geq \alpha$. Constraints of the first type are easy to deal with, since this is the same as to write $a_i+x_i\leq \alpha$ for all $i$, in terms of the usual arithmetic. Constraints of the second type mean that $a_i+x_i\geq \alpha$ for at least one $i$, and this can be written as 
$a_i+x_i+(1-w_i)M\geq \alpha$,
where $w_i\in\{0,1\}$ and $\sum_{i} w_i=1$, with $M$ a sufficiently large number. One can see that this reduction to MILP also applies to the constraints in~\eqref{e:newconstraints}.
Combining these techniques with the general Dempe-Franke algorithm is a matter of ongoing research.

Let us now discuss another approach to solving the problem
\begin{equation}
\label{e:DFrelaxed2}
\begin{split}
&\min_{x,y} f(x,y)\\
& \text{s.t.}\quad x\in\TPol_1,\quad y\in\TPol_2,
\quad x^Ty\leq \min\limits_{y'\in \TPol_2} x^Ty',
\end{split}
\end{equation}
where $f(x,y)$ is isotone with respect to the second argument:
$f(x,y^1)\leq f(x,y^2)$ whenever $y^1\leq y^2$. 
We can observe the following..


\begin{proposition}
\label{p:DFrelaxed}
If $f(x,y)$ is isotone with respect to the second argument then 
the minimum in~\eqref{e:DFrelaxed2} is equal 
to the minimum in the following problem:
\begin{equation*}
\begin{split}
&\min_{x,y} f(x,y)\\
& \text{s.t.}\quad x\in\TPol_1,\quad y\in\mathcal{M}(\TPol_{2}),
\quad x^Ty\leq \min\limits_{z\in \mathcal{M}(\TPol_2)} x^Tz.
\end{split}
\end{equation*}
\end{proposition}

This proposition provides for the following straightforward procedure 
solving~\eqref{e:DFrelaxed2} (and, in particular, Min-min Problem):

\begin{algorithm}
{\bf (Solving~\eqref{e:DFrelaxed2} and Min-min Problem)}\\[1.2 ex]
\label{a:DFStraight1}
{\rm 
\noindent {\bf Step 1.} Identify the set of 
minimal points $\mathcal{M}(\TPol_{2})$.\\[1.2 ex] 
{\bf Step 2.} For each point $y'\in\mathcal{M}(\TPol_{2})$
we solve the following optimization problem:
\begin{equation}
\label{e:DFStraight}
\begin{split}
&\min\limits_{x} f(x,y')\\
& \text{s.t.}\quad x\in\TPol_1, x^Ty'\leq x^Tz\quad\forall z\in \mathcal{M}(\TPol_2).
\end{split}
\end{equation}
{\bf Step 3.} Find the minimum among all 
problems~\eqref{e:DFStraight} for all $y'\in\TPol_2$.}$\hfill\square$
\end{algorithm}

Note that when $f(x,y)=a^Tx\oplus b^Ty$ for some vectors $a,b$ over $\Rmax$, 
problem~\eqref{e:DFStraight} can be solved by any algorithm of tropical linear programming~\cite{AllSimplex,But,GKS}. The set of all minimal points can be found by a combination of the tropical double description method of~\cite{AllDoubleD} that finds the set of all extreme points and the techniques of Preparata et al. for finding all minimal points of a finite set~\cite{PSh}, although clearly a more efficient procedure should be sought for this purpose.

\subsection{The max-max and min-max problems}

Let us now consider the problems where the lower-level objective 
is to maximize rather than to minimize:

\begin{equation}
\label{e:optmax}
\begin{split}
&\opt_{x,y} f(x,y)\\
& \text{s.t.}\quad x\in\TPol_1,\ y\in\arg\max_{y'}\{x^Ty'\colon y'\in\TPol_2\}.
\end{split}
\end{equation}

Following the LLVF approach, \eqref{e:optmax} is equivalent to
\begin{equation}
\label{e:optmaxLLVF}
\begin{split}
&\opt_{x,y} f(x,y)\\
& \text{s.t.}\quad x\in\TPol_1,\quad y\in\TPol_2,\quad x^Ty=\phi(x),
\end{split}
\end{equation}
where $\phi(x)=\max_z\{x^Tz\colon z\in\TPol_2\}$. The following
are similar to Definitions~\ref{d:minimal} and ~\ref{d:miness}.

\begin{definition}[Maximal Points]
\label{d:maximal}
Let $\TPol$ be a tropical polyhedron.
A point $x\in\TPol$ is called {\rm maximal} if 
$y\geq x$ and $y\in\TPol$ imply $y=x$.
\end{definition}

\begin{definition}[Max-Essential Subset]
\label{d:maxess}
Let $\TPol$ be a tropical polyhedron. Set $\mathcal{S}_{max}$ is called
a {\rm max-essential} subset of $\TPol$, if for any $x\in\Rmax^n$ the maximum
$\max_z\{x^Tz\colon z\in\TPol\}$ is attained at a point of $\mathcal{S}_{max}$.
\end{definition}

However, it is immediate that each compact tropical polyhedron contains 
its greatest point, and the above notions trivialize.

\begin{proposition}
\label{prop:unique-maximal}
Let $\TPol$ be a compact tropical polyhedron. Then $\TPol$ contains its greatest point $y^{\max}$.Furthermore, the singleton $\{y^{\max}\}$ is a max-essential subset of $\TPol$.
\end{proposition}

Proposition~\ref{prop:unique-maximal} implies that~\eqref{e:optmax}
(and~\eqref{e:optmaxLLVF}) are equivalent to 
\begin{equation}
\label{e:optmax-easy}
\begin{split}
&\opt_{x,y} f(x,y)\\
& \text{s.t.}\quad x\in\TPol_1,\quad y\in\TPol_2,\quad x^Ty=x^Ty^{\max},
\end{split}
\end{equation}
where $y^{\max}$ is the greatest point of $\TPol_2$.  The following result yields an 
immediate solution of the max-max problem.
\begin{corollary}[Solving Max-max Problem]
\label{c:maxmax-easy}
If $f(x,y)$ is isotone with respect to both arguments and 
$\opt=\max$, then $(x^{\max},y^{\max})$ is a globally optimal solution of~\eqref{e:optmax}, where $x^{\max}$ and $y^{\max}$ are the greatest points of $\TPol_1$ and $\TPol_2$.
\end{corollary}

Let us now consider~\eqref{e:optmax-easy} where $f$ is not necessarily isotone, 
or where $\opt=\min$ as in the case of Min-max problem. Suppose that 
$y^{\max}$ has all components in $\R$ and define point $x^*$ with coordinates:
\begin{equation*}
x_i^*=\bigotimes_{k\neq i} y_k^{\max}.
\end{equation*}
We first prove the following claim.

\begin{lemma}
\label{l:IJ}
Let $y^{\max}\in \R^n$. Consider sets $I$ and $J$ such that $I\cup J=[n]$ and $I\cap J=\emptyset$. Let $x$ be such that
\begin{equation}
\label{e:xIJ}
\begin{split}
x_i&=x_i^*\quad \forall i\in I,\\
x_i&<x_i^*\quad \forall i\in J.
\end{split}
\end{equation}
Then, if $y\in\TP_2$, equation $x^Ty=x^Ty^{\max}$ 
is equivalent to
\begin{equation}
\label{e:ymax2}
\bigoplus_{i\in I}\left(\bigotimes_{k\neq i}y_k^{\max}\right)  y_i
=\bigotimes_{k\in[n]} y_k^{\max}.
\end{equation}
\end{lemma}
\emph{Proof.}
Observe that $y\in\TP_2$ implies $y\leq y^{\max}$.
With such $x$ as in~\eqref{e:xIJ} and $y$ such that $y\leq y^{\max}$, we have
\begin{equation*}
\begin{split}
x^T y^{\max}&=\bigoplus_{i\in I} x_i^* y_i^{\max}\oplus \bigoplus_{j\in J} x_jy_j^{\max}=\bigoplus_{i\in I}\left(\bigotimes_{k\neq i}y_k^{\max}\right)  y_i^{\max}=\bigotimes_{k\in[n]} y_i^{\max},\\
x^T y&=\bigoplus_{i\in I} x_i^* y_i\oplus \bigoplus_{j\in J} x_jy_j
=\bigoplus_{i\in I} \left(\bigotimes_{k\neq i}y_k^{\max}\right)  y_i\oplus \bigoplus_{j\in J} x_jy_j.
\end{split}
\end{equation*}
Therefore, $x^T y=x^T y^{\max}$ becomes
\begin{equation}
\label{e:ymax1}
\bigoplus_{i\in I} \left(\bigotimes_{k\neq i}y_k^{\max}\right)  y_i\oplus \bigoplus_{j\in J} x_jy_j
=\bigotimes_{k\in[n]} y_k^{\max}.
\end{equation}
Moreover since $x_j<x_j^*$ we obtain that $x_jy_j< \bigotimes_{k\in[n]} y_k^{\max}(=x_j^*y_j^{\max})$ for 
each $j\in J$. Hence we can further simplify~\eqref{e:ymax1} to~\eqref{e:ymax2}.$\hfill\blacksquare$
 
Let us also introduce the following notation:

\begin{equation}
\label{e:TP1TP2IJ}
\begin{split}
\TPol_1^{IJ}&=\{x\in\TPol_1\colon  x_j(x^*_j)^{-1}<x_i(x_i^*)^{-1}\quad\forall i\in I, j\in J,\\
& x_k(x_k^*)^{-1}=x_l(x_l^*)^{-1}\quad\forall k,l\in I\}\\
\TPol_2^{IJ}&=\{y\in\TPol_2\colon \bigoplus_{i\in I} \left(\bigotimes_{k\neq i}y_k^{\max}\right)  y_i
=\bigotimes_{k\in[n]} y_k^{\max}\}
\end{split}
\end{equation}
Note that ``$x_j(x^*_j)^{-1}$'' means $x_j-x^*_j$ in the usual arithmetics.
Now, using Lemma~\ref{l:IJ} we can prove the following.

\begin{theorem}
\label{t:decomposition}
We have the following decomposition:
\begin{equation*}
\{(x,y)\in\TPol_1\times \TPol_2\colon x^Ty=x^Ty^{\max}\}
=\bigcup_{I,J} \{(x,y)\in \TPol_1^{IJ}\times \TPol_2^{IJ}\} 
\end{equation*}
where the union is taken over $I$ and $J$ are such that $I\cap J=\emptyset$ and $I\cup J=[n]$.
\end{theorem}

Theorem~\ref{t:decomposition} suggests that Problem~\eqref{e:optmax-easy} (and, equivalently,~\eqref{e:optmax}) can be solved by the following straightforward procedure. 

\begin{algorithm}
{\bf (Solving~\eqref{e:optmax} and Min-max Problem)}\\[1.2 ex]
\label{a:DFStraight2}
{\rm {\bf Step 1.} 
For each partition $I$, $J$ of $[n]$, identify the system of inequalities~\eqref{e:TP1TP2IJ} 
defining $\TPol_1^{IJ}$ and $\TPol_2^{IJ}$ and find a solution of the problem 
$\opt_{x,y} f(x,y)$ over 
$(x,y)\in \TPol_1^{IJ}\times\TPol_2^{IJ}$, if such solution exists.\\[1.2 ex] 
{\bf Step 2.} Compute $\opt$ over all solutions found at Step 1. $\hfill\square$}
\end{algorithm}

When $f(x,y)=a^Tx\oplus b^Ty$, this procedure reduces the problem to a finite number of tropical linear programming problems solved, e.g., by the algorithms of~\cite{AllSimplex,But,GKS}.

\begin{example}
\label{eg:2d}
{\rm Consider the following numerical example in two-dimensional case. 
Let $\TP_1$ is the tropical (max-plus) convex hull of the points $(-3,-1)$, $(-1,0)$ and $(-2,-3)$. See Figure \ref{fig:TP12} (a). $\TP_2$ is defined by $(1,1)$, $(0,0)$ and $(2,-1)$. See Figure \ref{fig:TP12} (b).

\begin{figure}
\centering
\label{fig:TP12}
\begin{tikzpicture}
\draw[blue] (-3,-1) -- (-2,-1);
\draw[green] (-2,-1) -- (-1,0);
\draw[red] (-2,-1) -- (-2,-3);
\draw[fill] (-3,-1) circle [radius=0.025];
\draw[fill] (-2,-1) circle [radius=0.025];
\draw[fill] (-2,-3) circle [radius=0.025];
\draw[fill] (-1,0) circle [radius=0.025];

\node [above] at (-1,0) {$(-1,0)$};
\node [left] at (-3,-1) {$(-3,-1)$};
\node [below] at (-2,-3) {$(-2,-3)$};
\node [right] at (-2,-1) {$(-2,-1)$};
\node at (-2,-4) {(a)};
\end{tikzpicture}
%
%
%
\begin{tikzpicture}
\coordinate (A) at (0,0);
\coordinate (B) at (1,1);
\coordinate (C) at (2,1);
\coordinate (D) at (2,0);
\coordinate (E) at (2,-1);

\draw (A) -- (B) -- (C) -- (D) -- (A);
\draw (D) -- (E);
\draw[red] (C) -- (E);
\draw[blue] (B) -- (C);
\draw[green]  (1,0.9) -- (1.9,0.9) -- (1.9,-1);

\draw[fill] (A) circle [radius=0.025];
\draw[fill] (B) circle [radius=0.025];
\draw[fill] (C) circle [radius=0.025];
\draw[fill] (D) circle [radius=0.025];
\draw[fill] (E) circle [radius=0.025];

\node [below] at (A) {$(0,0)$};
\node [left] at (B) {$(1,1)$};
\node [right] at (C) {$y^{max}$};
\node [below] at (E) {$(2,-1)$};
\node at (2,-2) {(b)};

\end{tikzpicture}
\caption{$\TP_1$ and $\TP_2$ of Example \ref{eg:2d}.}
\end{figure}


In this example, $y^{max} = (2,1)$ (the greatest point of $\TP_2$ in Figure \ref{fig:TP12} (b)). Therefore, $x^* = (1,2)$. Table \ref{tab:2d} shows three possible partitions of $\TP_1$ and $\TP_2$. Partition 1 corresponds to the line segment between $(-2,-1)$ and $(-2,-3)$ in $\TP_1$ and the line segment connecting $y^{max}$ and $(2,-1)$ in $\TP_2$ (red). Partition 2 corresponds to the line segment between $(-2,-1)$ and $(-3,-1)$ in $\TP_1$ and the line segment connecting $y^{max}$ and $(1,1)$ in $\TP_2$ (blue). Partition 3 corresponds to the line segment between $(-2,-1)$ and $(-1,0)$ in $\TP_1$ (green) and in $\TP_2$ the union of the line segment connecting $y^{max}$ and $(1,1)$ and the line segment between $y^{max}$ and $(2,-1)$ (green).}
\begin{table}
\label{tab:2d}
\caption{Partitions of Example \ref{eg:2d}}
\begin{tabular}{c c c c c}
    &$I$        &$J$         &$\TP_1^{IJ}$                                         &$\TP_2^{IJ}$        \\\hline
1   &\{1\}      &\{2\}       &\textcolor{red}{$\{x \in \TP_1: x_2-2<x_1-1\}$}           &\textcolor{red}{$\{y \in \TP_2\colon y_1 = 2\}$}                  \\
2   &\{2\}      &\{1\}       &\textcolor{blue}{$\{x \in \TP_1: x_1-1<x_2-2\}$}           &\textcolor{blue}{$\{y \in \TP_2\colon y_2 = 1\}$}                  \\
3   &\{1,2\}    &$\emptyset$ &\textcolor{green}{$\{x \in \TP_1: x_1-1=x_2-2\}$}           &\textcolor{green}{$\{y \in \TP_2\colon \max(1+y_1,\,2+y_2) = 3\}$}   
\end{tabular}
\end{table}

Assume the upper level objective is of the form $\min \; a^Tx \oplus b^Ty$, where $a$, $b \in \mathbb{R}^2$ . In ordinary algebra it can be written as $\min\;\{\max\{a_1+x_1,a_2+x_2,b_1+y_1,b_2+y_2\}\}$. It is obvious that the objective function is isotone with respect to $x$ and $y$. In partition 1, $x = (-2,-3)$ and $y = (2,-1)$ is always a solution regardless of $a$ and $b$. In partition 2, $x = (-3,-1)$ and $y = (1,1)$ is a solution. In partition 3, either $x = (-2,-1)$ and $y = (1,1)$ or $x = (-2,-1)$ and $y = (2,-1)$ solve the problem. However, these solutions are always dominated by the optimal points of partition 1 and partition 2. Therefore, in this example, it is sufficient to consider only partition 1 and partition 2. and decide between $(x,y)_1=((-2,-3),(2,-1))$
and $(x,y)_2=((-3,-1),(1,1)$. Taking $a_1=a_2=b_1=b_2$ makes $(x,y)_2$ an optimal solution of the problem, but taking $a_2=10$ and $a_1=b_1=b_2$ results in $(x,y)_1$.$\hfill\square$
\end{example}

\section{Conclusions}
We have studied the four different tropical analogues of a problem considered by Dempe and Franke~\cite{DF}. We showed that we can solve the problems by generalizing the Dempe-Franke algorithm and using reduction to MILP, or by decomposing the feasible set of a problem into a number of tropical polyhedra and performing tropical linear programming over these subdomains. The resulting methods need further practical study and theoretical improvement.

\section{Acknowledgement}
We gratefully acknowledge fruitful communication with Bart De Schutter and Ton van den Boom (TU Delft), who informed us about the reduction of tropical optimization problems to MILP.


\begin{thebibliography}{10}

\bibitem{AllSimplex}
Allamigeon, X., Benchimol, P., Gaubert, S., Joswig, M.:
Tropicalizing the simplex algorithm. SIAM J. on 
Discrete Math. \textbf{29}(2), 751-795 (2015).

\bibitem{AllDoubleD}
Allamigeon, X., Gaubert, S., Goubault, \'{E}.:
Computing the vertices of tropical polyhedra 
using directed hypergraphs. Discrete 
Comput. Geom. \textbf{49}, 247-279 (2013).

\bibitem{But}
Butkovi\v{c}, P.: Max-linear Systems: Theory and Algorithms. Springer, London (2010).

\bibitem{DF}
Dempe, S., Franke, S.:  Solution algorithm for an optimistic linear Stackelberg problem. Computers and Oper. Res. \textbf{41} 277-281 (2014).


\bibitem{Bart2} De Schutter, B., Heemels, W.P.M.H. and Bemporad, A.: On the equivalence of linear complementarity problems. Oper. Res. Letters, 
\textbf{30}(4) 211-222 (2002).

\bibitem{Bart3} De Schutter, B., Heemels, W.P.M.H., Bemporad, A.: Max-plus-algebraic problems and the extended linear complementarity problem -- algorithmic aspects.
In: Proceedings of the 15th IFAC World Congress, Barcelona, Spain (2002).

\bibitem{GKS}
Gaubert, S., Katz, R.D., Sergeev, S.: Tropical linear-fractional programming 
and parametric mean-payoff games. J. of Symb. Computation 
\textbf{47}(12) 1447-1478 (2012).

\bibitem{Helbig}
Helbig, S.: On Carath\'{e}odory's and Krein-Milman's theorems in fully ordered 
groups. Comment. Math. Univ. Carolin. \textbf{29}(1) 157-167 (1988).  


\bibitem{PSh} 
Preparata, F.P., Shamos, M.I.: Computational Geometry: An Introduction. 
Springer, New York (1985).


\end{thebibliography}
\end{document}